\documentclass[twoside,a4paper]{amsart}

\usepackage{courier} % tt
\normalfont
\usepackage{epsfig, graphicx, subfigure}
\usepackage{verbatim, setspace}
\usepackage{amsmath, amssymb}

\let\Re=\undefined
\DeclareMathOperator{\Re}{Re}
\let\Im=\undefined
\DeclareMathOperator{\Im}{Im}

\usepackage{amsmath,amssymb,latexsym,graphicx,amsthm}
\usepackage{latexsym}
\usepackage{eucal}
\usepackage{fullpage}
\usepackage{style}

\usepackage{tikz}
\newtheorem{theorem}{Theorem}[section]
\newtheorem{lemma}{Lemma}[section]

\theoremstyle{definition}

\newtheorem{open problem}{Open Problem}

\let\Re=\undefined
\DeclareMathOperator{\Re}{Re}
\let\Im=\undefined
\DeclareMathOperator{\Im}{Im}

\begin{document}
\title[Spatial asymptotics\ldots
 ]{Spatial asymptotics of Green's function and applications
  }

\author[S. Denisov]{Sergey A. Denisov}
\address{
\flushleft{
Department of Mathematics, University of Wisconsin-Madison,\\ 480 North Lincoln Dr., Madison, WI 53706, USA\\ 
and\\ 
Keldysh Institute of Applied Mathematics, Russian Academy of Science,\\ Miusskaya Pl. 4, Moscow, 125047, Russian Federation}}

\email{\href{mailto:denissov@wisc.edu}{denissov@wisc.edu}}

  \thanks{The work was supported by NSF DMS-1764245 and Van Vleck Professorship Research Award. 
 }

 \maketitle
\begin{abstract}

We study the spatial asymptotics of Green's function for the 1d  Schr\"odinger operator with  operator-valued decaying potential. The bounds on the entropy of the spectral measures are obtained. They are  used to establish the presence of a.c. spectrum.
\end{abstract} \vspace{1cm}

\section{Introduction and the main result}

In this note, we revisit the spectral theory of  Schr\"odinger operators with long-range potentials. In dimension one, the quest for the minimal assumptions on the decay of potential that guarantee the preservation of absolutely continuous spectrum resulted in the theorem (Deift-Killip \cite{dk}, see also \cite{ks1}), which says: 

{\it If $V\in L^2(\R^+)$, then
$\sigma_{\rm {ac}}(-\partial^2_{rr}+V)=[0,\infty)$ where $\sigma_{\rm {ac}}$ denotes a.c. spectrum of the operator with Dirichlet boundary condition at zero.}

In the case of the Dirac equation, an analogous result was obtained by M. Krein already in 1955 (see \cite{k1} and \cite{imrs}). $L^2$--condition is sharp: it is  known \cite{kls} that $V\in L^p(\R^+), p>2$ can lead to an empty a.c. spectrum. In higher dimension, one again is interested in finding the minimal assumptions on the decay of $V$ in $-\Delta+V, \,x\in \R^d, \,d\ge 2$ that guarantee ``scattering'' which can be understood either in the sense of preservation of a.c. spectrum or as  existence of wave operators in  Schr\"odinger dynamics. Some sharp results were obtained for decaying potentials that oscillate (see \cite{den1,saf1} for their surveys). However, if the oscillation assumption is dropped, then the identity $\sigma_{\rm ac}(-\Delta+V)=[0,\infty)$ is not known even for $V$ satisfying fairly strong conditions, e.g., $|V(x)|\le C(1+|x|)^{-1+\epsilon}, 0<\epsilon\ll 1$. Notice that the last assumption is only slightly weaker than the short-range condition of the classical scattering theory \cite{herm}. In this paper, we make progress on a related problem. \smallskip

Consider the Hilbert space $\mathcal{H}:=\oplus_{n=1}^\infty L^2(\R^+)$ with the inner product defined by 
\[
\langle F,G\rangle_{\mathcal{H}}=\int_0^\infty   \langle F,G\rangle dr=\sum_{n=1}^\infty \int_0^\infty f_n\overline{g_n}dr\,,
\]
where $F=(f_1,f_2,\ldots), G=(g_1,g_2,\ldots)$. 
We define the 1-d  Schr\"odinger operators
\begin{equation}\label{os}
H=-\partial^2_{rr}+V, \quad H^{(0)}=-\partial^2_{rr},\quad x \ge 0
\end{equation}
with Dirichlet boundary condition at the origin and operator-valued potential $V$. It satisfies $V(r)=V^*(r)$ for a.e. $r>0$ and $\|V\|\in L^\infty[0,\infty)$. By the general  theory of symmetric operators, $H$  defines the self-adjoint operator with the domain $\mathcal{D}(H)=\mathcal{D}(H^{(0)})=\oplus_{n=1}^\infty \mathcal{H}^2_0(\R^+)$, where $\mathcal{H}^2_0(\R^+):=\{f: f,f''\in L^2(\R^+), f(0)=0\}$ is the standard $\mathcal{H}^2(\R^+)$ Sobolev space of functions vanishing at the origin. Denote the Green's function of $H$ by $G(r,\rho,z)$, i.e.,
\[
R_zF=(H-z)^{-1}F=\int_{\R^+}G(r,\rho,z)F(\rho)d\rho\,, \,F\in \mathcal{H}\,.
\]
We let $z\in \C^+$ and $k=\sqrt z\in \{k\in \mathbb{C}^+,\Im k>0,\Re k>0\}$.
The Green's function of unperturbed operator will be called $G^{(0)}$.  Notice that
\begin{equation}\label{fi}
G^{(0)}(r,\rho,k^2)=\frac{i}{2k}\Bigl(e^{ik|r-\rho|}-e^{ik(r+\rho)}\Bigr)\,.
\end{equation}
Let $u:=R_{k^2}F,\, \psi:=e^{-ikr}u$.
We have
$
-u''+Vu=k^2u+F\,,\,  u(0,k)=0
$
and
\begin{equation}\label{l8}
-\psi''-2ik \psi'+V\psi=Fe^{-ikr}\,.
\end{equation}
In this note, we develop the perturbative theory which partially controls the spatial asymptotics of $u$ when $F$ has compact support, $r\to+\infty$ and $z\in \C^+$ is taken close to $\R^+$. Our analysis allows the direct study of $G(r,\rho,z)$ when $\rho$ is fixed and $r\to\infty$ but $u$ has a better regularity and we  will work with it instead. The following theorem showcases the typical  application of our analysis to the study of spectral type.

\begin{theorem}Suppose $\gamma>\frac 23, \lambda>0$, and $\|V\|\le \lambda(1+r)^{-\gamma}$. Then, $\R^+\subseteq \sigma_{\rm ac}(H)$.\label{maint}
\end{theorem}
Later in the text, we can assume that $\gamma$ is fixed in the range $\gamma\in (\frac 23,1)$. Many constants the reader encounters in this text  depend on $\gamma$ and $\lambda$ but we might not explicitly mention that.
\medskip

\noindent{\bf Remark. } The proof of the theorem employs elementary properties of subharmonic functions and a few  apriori integral estimates obtained directly from the equation itself. We avoid ODE asymptotical  methods so this technique can potentially be applied to study elliptic partial differential equations and difference operators on graphs. \smallskip

The connection between $\sigma_F$, the spectral measure of $F\in \mathcal{H}$, and the asymptotics of $u$ at infinity is revealed in the following lemma. 
\begin{lemma} Suppose $T>1$, $\supp V\subset  [0,T]$,\,  $F\in \mathcal{H}$, and $\supp F\subset [0,1]$. Then, $\sigma_F$ is absolutely continuous on $\R^+$ and
\begin{equation}\label{la3}
\sigma_F'(k^2)=k\pi^{-1}\|\psi(\infty,k)\|^2
\end{equation}
for $k\in \R^+$.
\end{lemma}
\begin{proof} Under the assumption of the lemma, the so-called absorption principle holds (see, e.g., \cite{gesz1,gesz2,gesz3} for the Weyl-Titchmarsh theory of operator-valued  Schr\"odinger operator).  In particular, for every interval $I\subset (0,\infty)$ and every positive $r$, the function $u(r,k)=(R_{k^2}F)(r)$
 has continuous extension in $k$ from $R_{I,1}:=I\times (0,1)$ to the interval $I$  and this $u$ satisfies $-u''+Vu=k^2u+F, \, u(0,k)=0$ for $k\in \overline{R_{I,1}}$. Thus, $\psi(r,k)=e^{-ikr}u(r,k)$ is defined as well for $k\in I$ and $\psi(r,k)=\psi(\infty,k)$ if $r>T$. That explains why the right-hand side in \eqref{la3} is well-defined. The absorption principle also implies that $\sigma_F$ is purely a.c. on $\R^+$.
Next, we take $k\in R_{I,1}$ and write
$
-u''+Vu=k^2u+F.
$
Take inner product with $u$ and integrate over $[0,T]$. Subtracting the resulting identity from its conjugate gives us
\[
\langle u'(T,k),u(T,k)\rangle-\langle u(T,k),u'(T,k)\rangle=(\bar k^2-k^2)\int_0^{T}\|u\|^2d\rho+\langle R_{k^2}F,F\rangle-\langle F,R_{k^2}F\rangle\,.
\]
Due to absorption principle, we can take $\Im k\to 0$ in the last formula. This gives \eqref{la3} after we take into account that $u(r,k)=e^{ikr}\psi(\infty,k)$ for $r>T$.

\end{proof}
\noindent {\bf Remark. }  One of the key ideas in the proof of Theorem \ref{maint} is based on the following observation. Taking the logarithm of the both sides in \eqref{la3} gives $\log \sigma_F'(k^2)=\log(k\pi^{-1})+2\log\|\psi(\infty,k)\|$. The function $\log\|\psi(\infty,k)\|$ is subharmonic in $R_{I,1}=I\times (0,1)$ for every closed interval $I\subset (0,\infty)$. Thus, rough bounds for $\log\|\psi(\infty,k)\|$ in $R_{I,1}$ can provide the lower bounds for the entropy $\int_{I'}\log \sigma'_F(k^2)dk, I'\subset I$ by application of mean-value inequality for subharmonic functions. The uniform control over the logarithmic integral implies the a.c. spectral type by the standard argument. A serious obstacle we will face is that the good  control of $\|\psi(\infty,k)\|$ is only possible when $\Im k$ is very small. The development of strategy that overcomes this difficulty was the main motivation to write this note. \smallskip

\noindent {\bf Some previous results.} In \cite{bs}, the reader can find an overview of one-dimensional results related to the topic. The survey papers \cite{den1,saf1} discuss the higher-dimensional case.  See also \cite{den2,saf2,saf3} for more recent advances. The one-dimensional  Schr\"odinger with operator-valued potential was extensively studied in the past and a thorough account of the literature can be found in \cite{gesz1,gesz2,gesz3}. The a.c. spectrum of operator-valued  Schr\"odinger with decaying potential was studied in the context of hyperbolic pencils in \cite{jfa}.

\smallskip
\noindent {\bf Motivation.} To relate \eqref{os} to multidimensional problems, consider the three-dimensional  Schr\"odinger operator $-\Delta+V, x\in \R^3$ which  allows the representation
\begin{equation}\label{kl1}
-\partial^2_{rr}-\frac{B}{r^2}+V(r,\theta)
\end{equation}
in the spherical coordinates $(r,\theta)\in \R^+\times \mathbb{S}^2$. Here, $B$ stands for Laplace-Beltrami operator on $\mathbb{S}^2$ and the Dirichlet boundary condition is assumed at the origin. If  the higher spherical modes can be neglected, one considers
\begin{equation}\label{kl2}
H=-\partial^2_{rr}-\frac{P_{\le r^\kappa}B}{r^2}+V(r,\theta)
\end{equation}
instead of \eqref{kl1}, where $P_{\le r^\kappa}$ is an orthogonal projection to the first $[r^\kappa]$ spherical harmonics. Assuming $|V(x)|\le C(1+|x|)^{-\gamma}$ with $\gamma>\frac 23$ and choosing $\kappa$ in a suitable way, we reduce \eqref{kl2} to the form \eqref{os}. 

\smallskip
\noindent {\bf Structure of the paper.} The second section  contains some apriori estimates for the solutions to equation \eqref{l8}. In the third section, we give the proof of  Theorem \ref{maint}. Some useful estimates on subharmonic functions are collected in Appendix 1. The second Appendix contains general bounds on Green's function.
\bigskip

\noindent {\bf Notation} \bigskip

$\bullet$ If $I$ is a closed interval on $\R$, $c_I$ denotes its center and $|I|$ denotes its length. $I_r$ stands for the interval centered at zero with radius $r$. $\R^+=(0,\infty)$.

$\bullet$ If $\psi$ is a vector in Hilbert space $\ell^2(\mathbb{N})$, then $\|\psi\|$ denotes its norm. If $V$ is a bounded linear operator acting in $\ell^2(\mathbb{N})$, then $\|V\|$ denotes its operator norm.

$\bullet$ If $I$ is an closed interval in $\R^+$ and $\delta>0$, then $R_{I,\delta}:=I\times (0,\delta)$.

$\bullet$ If $\phi,\psi\in \ell^2(\mathbb{N})$, then $\langle \phi,\psi\rangle $ refers to the inner product in $\ell^2(\mathbb{N})$.

$\bullet$ For $a>0$, we define $\log_+a=\max\{0,\log a\}, \,\log_-a=\min\{0,\log a\}$.

$\bullet$ The symbol $C_\alpha$ will indicate a positive constant whose dependence on a parameter  $\alpha$ we want to emphasize. The actual value of this constant can change from one formula to another.

$\bullet$  For two non-negative functions
$f_{1(2)}$, we write $f_1\lesssim f_2$ if  there is an absolute
constant $C$ such that
$
f_1\le Cf_2
$
for all values of the arguments of $f_{1(2)}$. We define $\gtrsim$
similarly and say that $f_1\sim f_2$ if $f_1\lesssim f_2$ and
$f_2\lesssim f_1$ simultaneously. If the constant $C$ depends on parameter $\alpha$, we might write $f_1\lesssim_\alpha f_2$.

$\bullet$ For the set $\Delta\subset \R$, we denote $\Delta^2=\{E^2: E\in \Delta\}$.
\bigskip

\section{Two simple  estimates obtained from the equation}

In this section, we consider the case when $\supp V\subset [0,T]$ and $\|V(r)\|<\lambda(r+1)^{-\gamma}, \gamma\in (\frac 23,1)$. In later discussion, we will be taking $T=2^n, n\ge n_0\gg 1$.
Let, e.g., $F$ be such that 
\begin{equation}\label{choice}
F=(f,0,\ldots), \quad  \|f\|_{L^2(\R^+)}=1,  \quad \supp f\subset [0,1], \quad f\nequiv 0\,.
\end{equation}
Let $\sigma_F$ be the spectral measure of $F$, i.e.,
\[
\langle R_zF,F\rangle_{\mathcal{H}}=\int \frac{d\sigma_F(E)}{E-z}\,,\, z\in \C\backslash \R\,.
\]
Recall that $\sigma_F$ is a probability measure and that $u=R_zF$. Rewrite equation \eqref{l8} for $\psi$ as
\begin{equation}\label{o11}
\psi'=i\frac{\psi''}{2k}-i\frac{V\psi}{2k}\,, \quad r>1\,.
\end{equation}
\begin{lemma}\label{lemmar} If $I$ is any closed interval in $\R^+,\, \alpha\in (0,1)$ and $k\in R_{I,T^{-\alpha}}$, then
\[
\sup_{r>0}\|\psi(r,k)\|\le C_{I,\alpha}\exp\left(2(\Im k)^{-(1-\alpha)/\alpha}\right)\,.
\]
\end{lemma}
\begin{proof} Since $V(r)=0$ for $r>T$, $\psi(r,k)=\psi(T,k)$ if $r>T$ and we can assume that $r\le T$. Because $\|u\|_{L^2[0,\infty)}\le C_I(\Im k)^{-1}$ we have $\|u''\|_{L^2[0,\infty)}\le C_I(\Im k)^{-1}$ from equation $-u''+Vu=k^2u+F$. Then, $\|u\|_{L^\infty[0,\infty)}\le C_I(\Im k)^{-1}$ as follows from the standard Sobolev's embedding.
%Lemma \ref{lemma11}.
Since $\|\psi(r,k)\|=e^{(\Im k) r}\|u(r,k)\|$, this gives us the statement of the lemma because
\[
(\Im k)r\le (\Im k)T\le (\Im k)^{-(1-\alpha)/\alpha}
\]
and 
\[
(\Im k)^{-1}\exp((\Im k)^{-(1-\alpha)/\alpha})\le C_{\alpha,I}\exp(2(\Im k)^{-(1-\alpha)/\alpha})
\]
\end{proof}
\noindent{\bf Remark.} Notice that this lemma only requires that $\|V\|\in L^\infty(\R^+)$ and $\supp V\subset [0,T]$.\bigskip

Next, we will study $\psi(r,k)$ when $r\in [T/2,T]$. In particular, we will be interested in how $\|\psi(r,k)\|$ deviates from $\|\psi(T/2,k)\|$ when $r>T/2$, $k\in R_{I,1}$, and $\Im k$ is small. Our basic tool is the following integral identity.

\begin{lemma} Let $1<a<b$ and $\Re k>0,\Im k> 0$, then
\begin{equation}\label{o2}
\|\psi(b,k)\|^2+\frac{\Im k}{|k|^2}\int_a^b \|\psi'\|^2d\rho=\|\psi(a,k)\|^2+Q_1-Q_2-\frac{\Im k}{|k|^2}\int_a^b \langle V\psi,\psi\rangle d\rho
\end{equation}
where
\[
Q_1:=\frac{i}{2k} \langle \psi'(b,k),\psi(b,k)\rangle -\frac{i}{2\bar k} \langle \psi(b,k),\psi'(b,k)\rangle 
\]
and 
\[
Q_2:=\frac{i}{2k} \langle \psi'(a,k),\psi(a,k)\rangle -\frac{i}{2\bar k} \langle \psi(a,k),\psi'(a,k)\rangle \,.
\]
\end{lemma}
\begin{proof}
Take inner product of both sides in \eqref{o11} with $\psi$ and integrate from $a$ to $b$. Then, take the real part of the resulting identity. We get
\[
\|\psi(b,k)\|^2=\|\psi(a,k)\|^2+\frac{i}{2k}\int_a^b \langle \psi'',\psi\rangle d\rho-\frac{i}{2\bar k}\int_a^b \langle \psi,\psi''\rangle d\rho-\frac{\Im k}{|k|^2}\int_a^b \langle V\psi,\psi\rangle d\rho\,.
\]
Integration by parts  gives the statement of the lemma.
\end{proof}
\noindent {\bf Remark.} In \eqref{o2}, an additional condition $V\ge 0$ immediately provides apriori estimate on $\int_1^\infty \|\psi'\|^2d\rho$ with essentially no  assumptions on the decay of $V$.\smallskip

The following lemma is straightforward.

\begin{lemma}\label{lemma1} Let $Y$ and $A$ be two $\ell^2(\mathbb{N})$-valued functions defined on $[a,\infty)$ that satisfy $\|Y\|,\|Y'\|,\|A\|\in L^2[a,\infty)$ and
\[
Y=\frac{i}{2k}Y'+A, \quad \Im k>0\,.
\]
Then, 
\begin{equation}\label{l1}
\|Y\|_{L^\infty[a,\infty)}\lesssim  \frac{|k|\|A\|_{L^2[a,\infty)}}{\sqrt{\Im k}}, \quad \|Y\|_{L^2[a,\infty)}\lesssim  \frac{|k|\|A\|_{L^2[a,\infty)}}{{\Im k}}.
\end{equation}
\end{lemma}
\begin{proof}
We have
$
Y'=-2ikY+2ikA 
$.  If $\Psi$ is defined by $\Psi:=e^{2ikr}Y$,
then
$
\Psi=-2ik \int_r^\infty A(s)e^{2iks}ds\,.
$
In the end, one has
\[
Y=-2ike^{-2ikr}\int_r^\infty A(s)e^{2kis}ds\,.
\]
Applying the convolution bounds, we get our lemma.
\end{proof}

If $T>1$, we arrange for two positive numbers $\mathcal{L}_
T$ and $\ell_T$ such that $\ell_T<\mathcal{L}_T$, $\ell_T:=T^{1-2\gamma+2\delta_1}$ and $\mathcal{L}_T:=T^{\gamma-1-\delta_1}$, where  $\delta_1$ is a positive parameter (e.g., take $\delta_1=\frac \gamma 2-\frac 13$). Its choice is possible since $\gamma\in (\frac 23,1)$. Given any closed interval $I\subset \R^+$, define the set
\[
PC_{I,T}:=R_{I,1}\cap \{k:  \ell_T\le \Im k \le \mathcal{L}_T\}\,.
\]
We will refer to  $PC_{I,T}$ as the {\it zone of perfect control}. The reader will see that this name is justified from the next two results.

\begin{lemma} For $k\in PC_{I,T/2}$, we have
\begin{eqnarray}\label{max2}
\|\psi(T,k)\|^2 =\|\psi(T/2,k)\|^2(1+\epsilon_T), \quad \epsilon_T\le C_I T^{-\delta_1}
 \end{eqnarray}
 where $\delta_1>0$.
\end{lemma}
\begin{proof}
We introduce 
$
M:=\sup_{r>T/2}\|\psi(r,k)\|\,.
$
Let $k\in R_{I,1}$. Applying Lemma~\ref{lemma1} to \eqref{o11} on the interval $[T/2,\infty)$, one has
\[
\|\psi'\|_{L^\infty [T/2,\infty)}\le C_I  M \frac{T^{0.5-\gamma}}{\sqrt{\Im k}}\,.
\]
Hence,
\[
\sup_{a,b>T/2}\|Q_{1(2)}\|\le C_I M^2\frac{T^{0.5-\gamma}}{\sqrt{\Im k}}\,.
\]
By the same Lemma \ref{lemma1}, 
\[
\|\psi'\|_{L^2[T/2,\infty)}\le C_I  M \frac{T^{0.5-\gamma}}{\Im k}\,.
\]
Taking supremum in $b\ge T/2$ in \eqref{o2} and letting $a=T/2$, we get
\[
|M^2- \|\psi(T/2,k)\|^2|\le C_I\Bigl((\Im k)T^{1-\gamma}+\frac{T^{0.5-\gamma}}{\sqrt{\Im k}}+\frac{T^{1-2\gamma}}{{\Im k}}\Bigr)M^2\,.
\]
Thus,
one has
\begin{equation}\label{max1}
M^2 =\|\psi(T/2,k)\|^2(1+\epsilon_T), \quad \epsilon_T\le C_I\Bigl((\Im k)T^{1-\gamma}+\frac{T^{0.5-\gamma}}{\sqrt{\Im k}}+\frac{T^{1-2\gamma}}{{\Im k}}\Bigr)\le C_I T^{-\delta_1}
\end{equation}
for given $k$. Now, we can take $b=T, a=T/2$ in \eqref{o2} and use the bound on $M$ to get the desired statement.
\end{proof}
We just saw that the $\|\psi(r,k)\|$ does not change much in $r$ when $r\in [T/2,T]$ and $k$ is fixed in {\it the zone of perfect control}. 
Next, we set up the iteration scheme which will play the key role in the proof of the main result.  Suppose $T_n=2^n, n\ge n_0$ where $n_0$ is a large parameter which will be fixed later. Given $V: \|V\|\le \lambda(1+r)^{-\gamma}, \gamma>\frac 23$, we let 
\begin{equation}\label{step1}
V_{(n)}:=V\cdot \chi_{[0,T_n]}, \quad H_{(n)}:=H^{(0)}+V_{(n)}, \quad \psi_n:=e^{-ikr}R_{(n),k^2}F\,,
\end{equation}
 where
function $F$ has been chosen in the beginning of this section and $R_{(n),z}:=(H_{(n)}-z)^{-1}$. The next lemma  estimates  $\psi_n(\infty,k)$ in the $(n-1)$-th {\it zone of perfect control}.

\begin{lemma} \label{lemmapc} Let $I$ be an closed interval in $\R^+$. If $k\in PC_{I,T_{n-1}}$, then
\[
\|\psi_{n}(\infty,k)\|=\|\psi_{n-1}(\infty,k)\|(1+\epsilon'_n),\quad  |\epsilon_n'|\le C_IT_n^{-\delta_2}
\]
where $\delta_2$ is a positive parameter.
\end{lemma}
\begin{proof} Recall that $\psi_j(T_j,k)=\psi_j(\infty,k)$ for every $j$.
By the previous lemma, it is enough to show that 
\begin{equation}\label{gy1}
\|\psi_{n}(T_n/2,k)\|=\|\psi_{n-1}(\infty,k)\|(1+O(T_n^{-\delta_3}))
\end{equation}
where $k\in PC_{I,T_{n-1}}$ and $\delta_3$ is a positive fixed number independent of $n$. To do that, we will use Lemma~\ref{l63}. Recall that $H_{(n)}=H_{(n-1)}+V\cdot \chi_{[T_{n-1},T_n]}$ and 
\[
R_{(n),k^2}F=R_{(n-1),k^2}F-R_{(n),k^2}  (V\cdot \chi_{[T_{n-1},T_n]}) R_{(n-1),k^2}F\,.
\]
Multiply the both sides with $e^{-ikr}$ and recall the definition of $\psi_n$ in \eqref{step1}. Since $\psi_{n-1}(r,k)=\psi_{n-1}(\infty,k)$ for $r\in [T_{n-1},\infty)$ and $k$ is in the {\it zone of perfect control}, we can apply Lemma \ref{l63} to $R_{(n),k^2}$. This yields
\begin{eqnarray}\label{gy2}
\|\psi_{n}(T_n/2,k)-\psi_{n-1}(\infty,k)\|\le C_I \|\psi_{n-1}(\infty,k)\|\int_{T_{n-1}}^{T_n} e^{(\Im k)(T_{n-1}-c(\rho-T_{n-1})-\rho)}T_n^{-\gamma}d\rho
\\\nonumber
\le C_I T_n^{-\gamma}(\Im k)^{-1}\|\psi_{n-1}(\infty,k)\|\le C_I T_n^{\gamma-1-2\delta_1}\|\psi_{n-1}(\infty,k)\|,
\end{eqnarray}
because $k\in PC_{I,T_{n-1}}$.  Putting together \eqref{gy1} and \eqref{gy2} gives the desired result.
\end{proof}
\bigskip

\section{Iteration and the proof of the main theorem}

Recall that $F$ is chosen to satisfy \eqref{choice}. First, we  need an auxiliary lemma.

\begin{lemma}Suppose $\|V\|\in L^\infty(\R^+)$ and $\psi_n$ is defined as in \eqref{step1}. Then,
\[
\sup_{0\le y\le 1}\int_I\|\psi_n(\infty,x+iy)\|^2dx<\infty, \quad \inf_{0\le y\le 1}\int_I\log\|\psi_n(\infty,x+iy)\|^2dx>-\infty
\]
for every closed interval $I\subset \R^+$.
\end{lemma}
\begin{proof}Since $V_{(n)}$ is compactly supported, $\psi_n(\infty,k)$ has continuous extension to any closed interval on the real line, and $\psi_n\nequiv 0$. It is also analytic in $k$ in every rectangle $R_{I,1}$ so the lemma follows from, e.g., the mean-value estimate for subharmonic function $\log\|\psi_n(\infty,k)\|$.
\end{proof}

To begin the iterative process which will be the key to the proof of our main result, we start with taking $I$, any closed interval in $\R^+$. Then, for this $I$, we choose
 $n_0\in \mathbb{N}$,  a fixed large parameter whose dependence on $I$ will be specified later, and define two numbers $A_{n_0}$ and $B_{n_0}$ as follows
 \begin{equation}\label{jumpstart}
A_{n_0}:=\sup_{0<y<\mathcal{L}_{T_{n_0}}}\int_{I}\|\psi_{n_0}(\infty,x+iy)\|^2dx, \quad B_{n_0}:=\sup_{0<y<\mathcal{L}_{T_{n_0}}}\int_{I}\log\|\psi_{n_0}(\infty,x+iy)\|dx\,.
 \end{equation}  
From the last lemma, one knows that $A_{n_0}<\infty$ and $B_{n_0}>-\infty$ for every $n_0$.  Next, we define the sequence of intervals $\{I_{(n)}\}, n\ge n_0$ by conditions
\begin{equation}\label{dok1}
I_{(n_0)}:=I,\, c_{I_{(n)}}=c_{I},\, |I_{(n)}|=|I_{(n-1)}|-2\tau_n
\end{equation}
and $\tau_n= T_{n}^{-\upsilon}$, where $0<\upsilon<0.01(-\gamma+1+\delta_1)$ so $\mathcal{L}_n=T_n^{\gamma-1-\delta_1}\ll \tau_n= T_{n}^{-\upsilon}$, see Figure 1. Notice that 
\[
\sum_{n\ge n_0}\tau_n=\sum_{n\ge n_0}2^{-\upsilon n}\sim C_\upsilon 2^{-\upsilon n_0}
\]
and $\lim_{n_0\to \infty}2^{-\upsilon n_0}= 0$. Therefore, if $I$ is given, we can always arrange for $n_0$ large enough that 
$
\mathcal{L}_{T_{n_0}}<1
$
and that there is $\widetilde I_{(n_0)}$:
\begin{equation}\label{sl}
c_{\widetilde I_{(n_0)}}=c_I,\quad  \widetilde I_{(n_0)}\subset \bigcap_{n\ge n_0}I_{(n)}, \quad \lim_{n_0\to \infty}|I\backslash \widetilde I_{(n_0)}|\to 0\,.
\end{equation}

\vspace{1cm}
\begin{tikzpicture}
\draw (0,1) -- (13,1);

\draw  [line width=0.35mm, red ]  (3.5,1) -- (9.5,1);

\draw (1,1) -- (1,4);
\draw (1,4) -- (12,4);
\draw (12,4) -- (12,1);
\draw (1,2) -- (12,2);

\draw[thick] (2,1) -- (2,3);
\draw[thick] (2,3) -- (11,3);
\draw[thick] (11,3) -- (11,1);
\draw[thick] (2,1.5) -- (11,1.5);
%\draw (2,1) -- (2,1.5);

%\draw (3,1) -- (3,1.5);
%\draw (2,1.5) -- (11,1.5);
%\draw (11,1.5) -- (11,1);

%\draw (10,1.5) -- (10,1);

\draw (0.4,2.4)   node[anchor=north]  {$\ell_{T_{n-1}}$};

\draw (12.6,4.2)   node[anchor=north]  {$\mathcal{L}_{T_{n-1}}$};

\draw (6,0.5)   node[anchor=north]  {Figure 1: $R_{I_{(n-1)},T_{n-1}}$ and $R_{I_{(n)},T_n}$};

%\draw (1.5,0.6)   node[anchor=north]  {$\delta$};

%\draw (2.5,0.6)   node[anchor=north]  {$\delta$};
\draw (6.5,1)   node[anchor=north]  {$c_{I_{(n-1)}}=c_{I_{(n)}}$};

\draw (8.5,1)   node[anchor=north]  {$\widetilde I_{(n_0)}$};

\draw (10.7,1)   node[anchor=north]  {$ I_{n}$};

\draw (12,1)   node[anchor=north]  {$ I_{(n-1)}$};

\draw[thick] (6.5,1)  circle (0.2mm);

\end{tikzpicture}
\bigskip

Let us collect  what we already know about the sequence $\{\psi_n\}$ below:\smallskip

$\bullet$ {\bf Rough upper bound}, Lemma \ref{lemmar}:
\begin{equation}\label{ap7}
\|\psi_n(\infty,k)\|\le C(I',\alpha)\exp\left(2(\Im k)^{-(1-\alpha)/\alpha}\right), \quad k\in R_{I',\mathcal{L}_{T_n}}\,,
\end{equation}
where $I'$ can be chosen as any open interval in $\R^+$ that contains $I_{(n_0)}=I$. The parameter $\alpha$ is related to $\gamma$ by $\alpha=1+\delta_1-\gamma$.\medskip

$\bullet$ {\bf The first step}: by construction, $A_{n_0}$ and $B_{n_0}$ are defined for every $n_0$.\medskip

$\bullet$ {\bf Estimate in the {\it zone of  perfect control}}, Lemma \ref{lemmapc}: if $k\in PC(I,T_{n-1})$, then
\begin{equation}\label{fact1}
\|\psi_{n}(\infty,k)\|=\|\psi_{n-1}(\infty,k)\|(1+\epsilon'_n),\quad  |\epsilon_n'|\le C_IT_n^{-\delta_2}\,.
\end{equation}\medskip

$\bullet$ {\bf Uniform bounds on the real line}, formula \eqref{la3}: for every $I'\subset \R^+$, we get

\begin{equation}\label{fact2}
\sup_{n\ge n_0}\int_{I'}\|\psi_{n}(\infty,k)\|^2dk<C_{I'}\,.
\end{equation}\medskip

To control $\psi_n(\infty,k)$  in $R_{I_{(n)},\mathcal{L}_{T_n}}$, one can use apriori estimates \eqref{ap7}, \eqref{fact2} along with \eqref{fact1}. To interpolate the  bounds on $\psi_n(\infty,k)$ from the {\it zone of perfect control} all the way to $R_{I_{(n)},\mathcal{L}_{T_n}}$, we will use a few estimates on the subharmonic functions that are collected and proved in the Appendix for reader's convenience. Our immediate goal is to prove the following lemma.
\begin{lemma}\label{fa} For every closed interval $J\subset \R^+$, we have the estimates
\begin{equation}\label{mest}
\limsup_{n\to\infty} \sup_{0<y<\mathcal{L}_{T_n}}\int_{J}\|\psi_n(\infty,x+iy)\|^2dx<\infty
\end{equation}
and 
\begin{equation}\label{mest4}
\|\psi_n(\infty,x+iy)\|^2\le   C_{J} \left(1+y^{-1}+(\mathcal{L}_{T_n}-y)^{-1}\right), \quad x\in J,\,\, 0<y<\mathcal{L}_{T_n}\,.
\end{equation}
\end{lemma}
\begin{proof}
We start with any interval $I$ and define the sequence $\{I_{(n)}\}$ as before in \eqref{dok1}. For each $n\ge n_0$, one lets
\[
A_{n}:=\sup_{0<y<\mathcal{L}_{T_n}}\int_{I_{(n)}}\|\psi_{n}(\infty,x+iy)\|^2dx\,.
\]
We will control how $A_n$  changes when $n$ is increased by one. 
Given $n-1$ and $A_{n-1}$, the goal is to estimate $A_{n}$.
To do that, we apply \eqref{fact1} and write
\begin{eqnarray*}
\sup_{\ell_{T_{n-1}}<y<\mathcal{L}_{T_{n-1}}}\int_{I_{(n-1)}}\|\psi_{n}(\infty,x+iy)\|^2dx\le (1+\epsilon'_n)^2\sup_{\ell_{T_{n-1}}<y<\mathcal{L}_{T_{n-1}}}\int_{I_{(n-1)}}\|\psi_{n-1}(\infty,x+iy)\|^2dx\\
\le A_{n-1}(1+\epsilon'_n)^2\,.
\end{eqnarray*}
Next, we apply \eqref{ap7}, \eqref{fact2}, and Lemma \ref{lems1} with $\kappa=(1-\alpha)/\alpha, \delta\sim \tau_{n}, \epsilon_1\sim \mathcal{L}_{T_{n-1}}, \epsilon_2=\ell_{T_{n-1}}$ 
to get
\begin{eqnarray*}
\sup_{0<y<\ell_{T_{n-1}}}\int_{I_{(n)}}\|\psi_{n}(\infty,x+iy)\|^2dx\le  C_{I'}+O(T_n^{-\delta_4}(1+C_{I'}+A_{n-1})), \quad \delta_4>0\,.
\end{eqnarray*}
In the end, we have
\[
A_{n}\le \max\Bigl\{C_{I'}+O(T_n^{-\delta_4}(1+C_{I'}+A_{n-1})), A_{n-1}(1+O(T_n^{\delta_5}))\Bigr\}
\]
with positive $\delta_4$ and $\delta_5$. That is supplemented by fixing $A_{n_0}$. The previous bound yields
\[
A_{n}\le A_{n-1}(1+O(T_n^{\delta_5}))+O(T_n^{-\delta_4})
\]
and $A_n\le C_I A_{n_0}$\,. Consequently, 
\begin{equation}\label{mesto}
\limsup_{n\to\infty} \sup_{0<y<\mathcal{L}_{T_n}}\int_{\widetilde I_{(n_0)}}\|\psi_n(\infty,x+iy)\|^2dx<\infty\,.
\end{equation}
Due to \eqref{sl}, we can start with any $J$, choose $I$ that contains it and then $n_0$ so large that $\widetilde{I}_{(n_0)}$ contains $J$ too. That  will give us the first statement of the lemma. Now, the bound \eqref{mest4} follows  from \eqref{k2}.
\end{proof}

\begin{lemma} For every closed interval $J\subset \R^+$, we have an estimate
\begin{equation}\label{mest1}
\liminf_{n\to\infty} \int_{J}\log \|\psi_n(\infty,k)\|dk>-\infty\,.
\end{equation}
\end{lemma}

\begin{proof}As in the previous proof, we define
\[
\quad B_{n}:=\inf_{0<y<\mathcal{L}_{T_n}}\int_{I_{(n)}}\log \|\psi_{n}(\infty,x+iy)\|dx\,.
\]
We will control how  $B_n$ changes when $n$ is increased by one. 
Given $B_{n-1}$ and the previous lemma, we want to estimate $B_{n}$. To control $\log \|\psi_{n}(\infty,x+iy)\|$ 
in the upper part of $R_{I_{(n)},T_n}$, we use estimates in
 $PC(I_{n-1},T_{n-1})$.
Applying  \eqref{fact1}, one has
\begin{eqnarray}\label{int1}
\inf_{\ell_{T_{n-1}}\le y\le \mathcal{L}_{T_{n}}}\int_{I_{(n-1)}}\log\|\psi_{n}(\infty,x+iy)\|dx=\\ O(\epsilon'_n)+\inf_{\ell_{T_{n-1}}\le y\le \mathcal{L}_{T_{n}}}\int_{I_{(n-1)}}\log \|\psi_{n-1}(\infty,x+iy)\|dx\ge O(\epsilon'_n)+B_{n-1}.\nonumber
\end{eqnarray}
and
\begin{eqnarray*}
\inf_{\ell_{T_{n-1}}\le y\le \mathcal{L}_{T_{n}}}\int_{I_{(n)}}\log\|\psi_{n}(\infty,x+iy)\|dx= O(\epsilon'_n)+\inf_{\ell_{T_{n-1}}\le y\le \mathcal{L}_{T_{n}}}\int_{I_{(n)}}\log \|\psi_{n-1}(\infty,x+iy)\|dx\,.
\end{eqnarray*}
Notice that for the chosen range of $y$ we have
\[
\int_{I_{(n)}}\log \|\psi_{n-1}(\infty,x+iy)\|dx= \int_{I_{(n-1)}}\log \|\psi_{n-1}(\infty,x+iy)\|dx-\int_{I_{(n-1)}\backslash I_{(n)}}\log\|\psi_{n-1}(\infty,x+iy)\|dx
\]
and
\[
-\int_{I_{(n-1)}\backslash I_{(n)}}\log\|\psi_{n-1}(\infty,x+iy)\|dx\ge -\int_{I_{(n-1)}\backslash I_{(n)}}\log_+ \|\psi_{n-1}(\infty,x+iy)\|dx.
\]
Then, 
\[
\int_{I_{(n-1)}\backslash I_{(n)}}\log_+ \|\psi_{n-1}(\infty,x+iy)\|dx\lesssim_I  \tau_n^{\frac 12 }
\]
as follows from the estimate $\log_+t\le |t|$, Cauchy-Schwarz inequality, \eqref{mest}, and the bound $|I_{(n)}\backslash I_{(n-1)}|\lesssim \tau_n$. In the end, we get
\[
\inf_{\ell_{T_{n-1}}\le y\le \mathcal{L}_{T_{n}}}\int_{I_{(n)}}\log\|\psi_{n}(\infty,x+iy)\|dx\ge B_{n-1}+O(\tau_n^{\frac 12})+O(\epsilon'_n)\,.
\]
To control the integral for the smaller values of $y$, i.e., when $y<\ell_{T_{n-1}}$, we apply  Lemma \ref{lems2} with $\epsilon_1=\mathcal{L}_{T_n}, \epsilon_2=2\ell_{T_{n-1}}$ and $\delta\sim \tau_n$. The base of smaller rectangle is $I_{(n)}$ and the base of the larger one is $I_{(n-1)}$. Given Lemma \ref{fa}, we can write
\begin{eqnarray*}
\inf_{0<y<\ell_{T_{n-1}}}\int_{I_{(n)}}\log \|\psi_{n}(\infty,x+iy)\|dx\ge (1+O(T_n^{-\delta_6}))\int_{I_{(n-1)}}\log \|\psi_{n}(\infty,x+2i\ell_{T_{n-1}})\|dx-O(T_n^{-\delta_7})\,.
\end{eqnarray*}
with positive $\delta_6$ and $\delta_7$. For the integral on the right-hand side,  apply \eqref{int1}.  In the end, one has
\[
B_{n}\ge (1+O(T_n^{-\delta_8}))B_{n-1}+O(T_n^{-\delta_9})\,, \quad \delta_8>0,\delta_9>0.
\]
Consequently $\liminf_{n\to\infty} B_n>-\infty$ and thus
$
\liminf_{n\to\infty}\int_{I_{(n)}}\log\|\psi_n(\infty,x)\|dx>-\infty.
$
Since 
\[
\int_{I_{(n)}}\log\|\psi_n(\infty,x)\|dx=\int_{I_{(n)}}\log_-\|\psi_n(\infty,x)\|dx+\int_{I_{(n)}}\log_+\|\psi_n(\infty,x)\|dx
\]
and \eqref{fact2} guarantees that
$
\limsup_{n\to\infty}\int_{I_{(n)}}\log_+\|\psi_n(\infty,x)\|dx<\infty\,,
$
we have
\[
\liminf_{n\to\infty}\int_{\widetilde I_{(n_0)}}\log_-\|\psi_n(\infty,x)\|dx\ge
\liminf_{n\to\infty}\int_{I_{(n)}}\log_-\|\psi_n(\infty,k)\|dk>-\infty.
\]
The reasoning given at the end of the proof of the previous lemma can be used again to deduce \eqref{mest1}.
\end{proof}

The last two results provide the crucial estimates for $\|\psi_n(\infty,k)\|$ when $\Im k\in (0,\mathcal{L}_{T_n})$. They control the behavior of $\|(R_{(n),k^2}F)(r)\|$ for large $r$ without giving precise asymptotics for $(R_{(n),k^2}F)(r)$. That, however, is enough to prove Theorem \ref{maint}.\smallskip

{\it Proof of Theorem \ref{maint}.}
Take any closed interval $J\subset \R^+$  and recall that  $V_{(n)}=V\cdot \chi_{r<T_n}$. Define $\sigma_{(n),F}$, the spectral measure of $F$ relative to $H_{(n)}=H^{(0)}+V_{(n)}$. The spectral measure of $F$ relative to $H$ is  $\sigma_F$. Then, the previous lemma yields
\[
\liminf_{n\to\infty}\int_{{\Delta}^2}\log \sigma_{(n),F}'(E)dE>-\infty\,.
\]
Since $\lim_{n\to\infty}\|R_{(n),z}F-R_{z}F\|_{\mathcal{H}}=0, \,\, z\in \mathbb{C}^+$, we get $\sigma_{(n),F}\to \sigma_F$ in the weak--$(\ast)$ sense. Hence, (see \cite{ks}, section~5),
\[
\int_{\Delta^2}\log \sigma'_FdE>-\infty
\]
which implies that $\Delta^2$ supports a.c. spectrum of the original $H$. Since $\Delta$ was arbitrary, we get the statement of the theorem. \qed
\bigskip

\section{Appendix 1: some estimates on subharmonic functions}

For the reader's convenience, we collect  some elementary estimates on subharmonic functions in this appendix. Start with the estimates for the subharmonic function of a thin isosceles trapezoid. We denote this trapezoid by ${\mathcal{T}}_{I,\epsilon,\beta}$  where the height is $\epsilon$,  the side angles at the lower base are both equal to $\pi/\beta$, and the projection of the upper base to the real line is a given interval $I\subset \R$.
First, we will need some estimates on the harmonic measure of that trapezoid. It is instructive to start with giving the exact formula for harmonic measure of the infinite tube which is ``infinitely long'' rectangle. If $Cyl_{\epsilon}:=\{k: 0<\Im k<\epsilon\}$, then the density of harmonic measure on its lower side is
\begin{equation}\label{lk1}
\omega'_{k}(t)=\frac{1}{2\epsilon} \frac{\sin (\pi \epsilon^{-1}y)}{\cosh(\pi \epsilon^{-1}(x-t))-\cos (\pi \epsilon^{-1}y)}\,, \, t\in \R, \, k=x+iy\in Cyl_{\epsilon}\,.
\end{equation}
That formula can be verified directly.  Let $\Gamma:=\partial{\mathcal{T}_{J,\epsilon,\beta}}=\Gamma_1\cup\ldots\cup\Gamma_4$, where $\Gamma_1$ is an upper base, $\Gamma_2$ the lower base, $\Gamma_3$ the left leg, and $\Gamma_4$ the right leg of the trapezoid. 
Denote the harmonic measure at point $k$ by $\omega_k$.

\begin{lemma} Suppose the $\Gamma_2=[0,2]$ and the positive parameters $\beta,\epsilon,\delta$ are chosen such that $\beta>2, \beta\sim 1\,,\, \epsilon< \delta^2\ll 1$, $k=x+iy \in R_{(\delta,2-\delta), 0.5\epsilon}$, and $\xi\in \Gamma$. Then, the derivative of harmonic measure in the corresponding trapezoid with respect to its arclength satisfies
\begin{eqnarray}
\xi=s+\epsilon i\in \Gamma_1, \hspace{3cm}
\label{hm1}
\omega'_k(\xi)\lesssim    \frac{\epsilon^{-2}y}{\cosh (\pi \epsilon^{-1}(x-s))} \,,\\
\xi=s\in \Gamma_2, \quad\hspace{3cm}
\label{hm2}
\omega'_k(\xi)\le   \frac{y}{\pi((s-x)^2+y^2)}\,,
\\
\xi=te^{i\pi/\beta}\in \Gamma_3, \hspace{3cm}\hspace{1cm}\quad \label{hm3}
\omega'_k(\xi)\le C_{\beta}  \frac{(xt)^{\beta-1}y}{(t^2+x^2)^\beta}\,,\\
\xi=2+te^{i(\pi-\pi/\beta)}\in \Gamma_4,\, x<1\quad \quad\hspace{3cm} \label{hm4}
\omega'_k(\xi)\le C_{\beta} yt^{\beta-1}\,.
\end{eqnarray}
\end{lemma}
\begin{proof} See Figure 2.

\begin{tikzpicture}
\draw (0,1) -- (15,1);

%\draw (1.5,1) -- (3,4);
%\draw (13.5,1) -- (12,4);

\draw (0.5,1) -- (3,4);
\draw (14.5,1) -- (12,4);

\draw (3,1) -- (3,4);

\draw (4,1) -- (4,4);

\draw (3,4) -- (12,4);
\draw (12,4) -- (12,1);

\draw (11,1) -- (11,4);

%\draw (3,1) -- (4,3);

%\draw (3,2) -- (12,2);
\draw (12,3) -- (12,1);

\draw (7,0.5)   node[anchor=north]  {Figure 2, $\epsilon< \delta^2\ll 1$};

\draw (5.1,1.5)   node[anchor=north]  {$k$};
\draw[thick] (5,1.5)  circle (0.2mm);

\draw[thick] (7.5,1)  circle (0.2mm);
\draw (7.5,1)   node[anchor=north]  {$1$};

%\draw (2.2,1.5)   node[anchor=north]  {$\pi/(2\beta)$};
%\draw (12.8,1.5)   node[anchor=north]  {$\pi/(2\beta)$};

\draw (0.7,1)   node[anchor=north]  {$0$};
\draw (14.3,1)   node[anchor=north]  {$2$};
\draw (4,1)   node[anchor=north]  {$\delta$};
\draw (12,1)   node[anchor=north]  {$2-\delta$};
\draw (1.2,1.5)   node[anchor=north]  {$\pi/\beta$};
\draw (13.8,1.5)   node[anchor=north]  {$\pi/\beta$};

%\draw (11.7,1.7)   node[anchor=north]  {$\epsilon_2$};

\draw (3.3,2.6)   node[anchor=north]  {$\epsilon$};

%\draw (11,1)   node[anchor=north]  {$1$};

%\draw (12,1)   node[anchor=north]  {$1+\delta$};

\draw (7.5,4.7)   node[anchor=north]  {$\Gamma_1$};

\draw (7.5,1.7)   node[anchor=north]  {$\Gamma_2$};

\draw (1.5,3.2)   node[anchor=north]  {$\Gamma_3$};

\draw (13.5,3.2)   node[anchor=north]  {$\Gamma_4$};

\end{tikzpicture}
\vspace{0.5cm}

Recall the following monotonicity property of harmonic measure. If $\Omega_1\subset \Omega_2$ and $E\subset \partial\Omega_1\cap \partial\Omega_2$, then $\omega_{k,\Omega_1}(E)\le \omega_{k,\Omega_2}(E)$ for $k\in \Omega_1$ (\cite{Garnett}, p. 36) where $\omega_{k,\Omega}$ denotes harmonic measure at point $k$ relative to the domain $\Omega$. This monotonicity helps us get the required upper bounds by comparing to harmonic measure of an angle, an infinite cylinder, or a half-plane. We obtain \eqref{hm1}  by comparing with infinite cylinder and \eqref{hm2} by comparing with the upper half-plane. The other two formulas are deduced by making a comparison with an infinite angle.
\end{proof} 

\noindent {\bf Remark.} The estimates in the upper part of rectangle can be obtained in a similar way.\smallskip

We will  need the following result later. Recall that $I_r$ denotes the interval on the real line with radius $r$ centered at the origin.
\begin{lemma}\label{new_l1} Suppose the positive parameters $\epsilon_2,\epsilon_1,\delta$ satisfy $2\epsilon_2<\epsilon_1< \delta^{2}\ll 1$ and let $\omega_k$ be a harmonic measure for $R_{I_{1+\delta},\epsilon_1}$. Then, for $k=x+i\epsilon_2$, we have 
\begin{equation}\label{new1}
\sup_{|\xi|<1-\delta}\left|\int_{I_{1}}\omega_{x+i\epsilon_2}'(\xi)dx-1\right|\lesssim \epsilon_2\epsilon_1^{-1}\,.
\end{equation}
\end{lemma}
\begin{proof} The required density of harmonic measure can be written via harmonic measure of infinite cylinder through proper extension from $I_{1+\delta}$ to $\R$. The resulting formula shows that the contribution from the left and righ sides of rectangle are exponentially small and the desired density can be well approximated by the density of harmonic measure of the infinite cylinder. Then, we use formula \eqref{lk1} to obtain required bound.
\end{proof}

\begin{lemma} \label{lems1} Suppose the positive parameters $\epsilon_1,\epsilon_2$ and $\delta$ satisfy $2\epsilon_2<\epsilon_1< \delta^{2}\ll 1$.
 Assume that $h$ is $\ell^2(\mathbb{N})$-valued function  holomorphic in  $R_{I_{2},1}$, continuous in  $\overline{R_{I_{2},1}}$, and 
\begin{equation}\label{k1}
\|h(k)\|\le  C_1\exp(C_2(\Im k)^{-\kappa}), \,\, k\in R_{I_{2},1}\,,\,  1<\kappa,\, \kappa\sim 1.
\end{equation}
Then,  we have
\begin{eqnarray}\label{k2}
\|h(x+iy)\|^2\le   C_\kappa  \left(1+y^{-1}A+(\epsilon_1-y)^{-1}B\right),\\ A:=\int_{ I_{2} }\|h(t)\|^2dt, \quad B:=\int_{-1-\delta }^{1+\delta}\|h(t+i\epsilon_1)\|^2dt\nonumber
\end{eqnarray}
provided that $k=x+iy\in R_{I_{1+\frac \delta 2 },\epsilon_1}$.
Moreover,
\begin{eqnarray}\label{k3}
\sup_{0<y<\epsilon_2} \int_{-1}^1\|h(x+iy)\|^2dx\le A
+C_\kappa \epsilon_2\epsilon_1^{-1}(A+B+\epsilon_1)\,.
\end{eqnarray}

\end{lemma}
\begin{proof} See Figure 3.

\begin{tikzpicture}
\draw (0,1) -- (15,1);

%\draw (1.5,1) -- (3,4);
%\draw (13.5,1) -- (12,4);

\draw (0.5,1) -- (3,4);
\draw (14.5,1) -- (12,4);

\draw (3,1) -- (3,4);

\draw (4,1) -- (4,4);

\draw (3,4) -- (12,4);
\draw (12,4) -- (12,1);

\draw (11,1) -- (11,4);

%\draw (3,1) -- (4,3);

\draw (3,2) -- (12,2);
\draw (12,3) -- (12,1);

\draw (7,0.5)   node[anchor=north]  {Figure 3};

\draw (5.1,1.5)   node[anchor=north]  {$k$};
\draw[thick] (5,1.5)  circle (0.2mm);

\draw[thick] (7.5,1)  circle (0.2mm);

\draw (7.5,1)   node[anchor=north]  {$0$};

%\draw (2.2,1.5)   node[anchor=north]  {$\pi/(2\beta)$};
%\draw (12.8,1.5)   node[anchor=north]  {$\pi/(2\beta)$};

\draw (0.7,1)   node[anchor=north]  {$\pi/(2\kappa)$};
\draw (14.3,1)   node[anchor=north]  {$\pi/(2\kappa)$};

\draw (11.7,1.7)   node[anchor=north]  {$\epsilon_2$};

\draw (3.3,2.6)   node[anchor=north]  {$\epsilon_1$};

\draw (11,1)   node[anchor=north]  {$1$};

\draw (12,1)   node[anchor=north]  {$1+\delta$};

\draw (4,1)   node[anchor=north]  {$-1$};

\draw (3,1)   node[anchor=north]  {$-(1+\delta)$};

\draw (7.5,4.7)   node[anchor=north]  {$\Gamma^+$};
\end{tikzpicture}
\vspace{0.5cm}

We can assume $h\nequiv 0$.
 Let $k=x+iy\in R_{I_{1+\frac \delta 2 },\epsilon_1}$. Consider the isosceles trapezoid ${\mathcal{T}}_{{I_{1+\delta}},\epsilon_1,\pi/(2\kappa)}$. Denote its upper base by
 $\Gamma^+$ and its lower base by $\Gamma^-$. 
We write the mean-value inequality for subharmonic function $2\log_+\|h\|$ and use the estimate \eqref{hm3} on the density of harmonic measure on the legs to get
\begin{eqnarray*}
2\log_+\|h(k)\|\le 2\int_{\partial {\mathcal{T}}_{{I_{1+\delta},\epsilon_1},\pi/(2\kappa)}}\log_+\|h\|d\omega_k\le C_\kappa y\delta^{-1-2\kappa}\epsilon_1^{\kappa}+\\
2\int_{\Gamma^+\cup \Gamma^-}\log_+ \|h\|\omega_k'(\xi)d\xi\le C_\kappa +2\int_{\Gamma^+\cup\Gamma^-}\log_+ \|h\|\omega_k'(\xi)d\xi
\end{eqnarray*}
where we applied the given estimates on $\|h\|$ along with   $\epsilon_1<\delta^{2}$. Define $Q(k)=\max\{1,\|h\|\}$ and notice that $\log Q=\log_+ Q\ge 0$ so
\[
\log Q^2\le C_\kappa+\int_{\Gamma^+\cup\Gamma^-}(\log Q^2)\omega_k'(\xi)d\xi\le C_\kappa+\int_{\Gamma^+\cup\Gamma^-}(\log Q^2)d\mu, \quad\mu:=\frac{\omega_k|_{\Gamma^-\cup \Gamma^+}}{\|\omega_k|_{\Gamma^-\cup \Gamma^+}\|}\ge \omega_k|_{\Gamma^-\cup \Gamma^+}\,.
\]
Taking the exponential of both sides and using Jensen's inequality
\[
\exp \left(\int \log f {d\mu}\right)\le \int f{d\mu}, \quad \|\mu\|=1
\]
we get
\[
Q^2\le C_\kappa  \frac{\displaystyle \int_{ \Gamma^-\cup  \Gamma^+  } Q^2\omega_k'(\xi)d\xi}{\|\omega_k|_{\Gamma^-\cup \Gamma^+}\|}\,.
\]
For considered $k$, we have $\|\omega_k|_{\Gamma^-\cup \Gamma^+}\|\sim 1$. Thus,
\begin{eqnarray*}
Q^2
\lesssim_\kappa  1+\int_{-2}^2 \frac{\pi^{-1}y}{(\xi-x)^2+y^2}\|h(\xi)\|^2d\xi 
+\int_{-1-\delta}^{1+\delta} \frac{\pi^{-1}(\epsilon_1-y)}{(\xi-x)^2+(\epsilon_1-y)^2}\|h(\xi+i\epsilon_1)\|^2d\xi 
\\
\lesssim_\kappa 1+C\left(y^{-1}\int_{I_2}\|h\|^2d\xi+(\epsilon_1-y)^{-1}\int_{-1-\delta}^{1+\delta} \|h(\xi+i\epsilon_1)\|^2d\xi \right)\,.
\end{eqnarray*}
To obtain \eqref{k3}, we take $k\in R_{I_1,\epsilon_2}$ and apply the mean-value inequality to subharmonic function $\|h(k)\|^2$ inside the domain  $R_{I_{1+\frac \delta 2 },\epsilon_1}$. The symbol 
$\Gamma^+_{I_{1+\frac \delta 2 },\epsilon_1}$ will stand for an upper base of this rectangle. Then,
\begin{equation}\label{k21}
\|h(k)\|^2\le \int_{\partial R_{I_{1+\frac \delta 2 },\epsilon_1}}\|h\|^2d\omega_k\le I+\int_{I_{1+\frac \delta 2 }} \|h\|^2\omega_k'(\xi)d\xi+\int_{\Gamma^+_{I_{1+\frac \delta 2 },\epsilon_1}}\|h\|^2\omega_k'(\xi)d\xi\,.
\end{equation}
To estimate the first term, we use \eqref{k2}. That gives
\begin{eqnarray*}
I\lesssim_\kappa \int_0^{0.5 \epsilon_1 }   \left(1+t^{-1}A+(\epsilon_1-t)^{-1}B\right)
\left( \frac{xyt}{(x^2+t^2)^2} \right)dt+ \\
  \int_{0.5 \epsilon_1 }^{\epsilon_1}
\left(1+t^{-1}A+(\epsilon_1-t)^{-1}B\right)
\left( \frac{xy(\epsilon_1-t)}{(x^2+(\epsilon_1-t)^2)^2} \right)dt
\lesssim_\kappa    (A+B+\epsilon_1)y\epsilon_1 \delta^{-3}
\end{eqnarray*}
as follows from \eqref{k2} and the estimates for the harmonic measure of rectangle. For the last term in the right hand side of \eqref{k21}, one employs the bound on harmonic measure to write
\begin{equation}\label{gt1}
\int_{\Gamma^+_{I_{1+\frac \delta 2 },\epsilon_1}}
\|h\|^2\omega_k'(\xi)d\xi\lesssim\int_{-1-\frac \delta 2 }^{1+\frac \delta 2 } \|h(\xi+i\epsilon_1)\|^2  \frac{\epsilon_1^{-2}y}{\cosh (\pi \epsilon_1^{-1}(x-\xi))}  d\xi\,.
\end{equation}
Next, we integrate  \eqref{k21} in $x\in I_1$. 
Integration of \eqref{gt1} yields
\[
\int_{I_1}\left(\int_{\Gamma^+_{I_{1+\frac \delta 2 },\epsilon_1}}
\|h\|^2\omega_k'(\xi)d\xi\right)dx\lesssim \int_{-1-\frac \delta 2 }^{1+\frac \delta 2 } \|h(\xi+i\epsilon_1)\|^2\left(\int_{I_1} \frac{\epsilon_1^{-2}y}{\cosh (\pi \epsilon_1^{-1}(x-\xi))}  dx\right)d\xi\lesssim By\epsilon_1^{-1}\,.
\]
The second term on the right-hand side of \eqref{k21}  contributes
\[
\int_{I_1}\left(\int_{I_{1+\frac \delta 2 }} \|h\|^2\omega_k'(\xi)d\xi\right)dx
\le \int_{I_{1+\frac \delta 2 }}\|h\|^2\left(\int_{I_1}\omega_k'(\xi) d\xi \right)dx\le  \int_{ I_2}\|h\|^2dx
\]
where the estimate
\[
\omega_k'(\xi)\le \frac{\pi^{-1}y}{(\xi-x)^2+y^2}
\]
was used.
Combining the bounds, we get \eqref{k3} after our assumption $\epsilon_1<\delta^2$ is taken into account.
\end{proof}

\begin{lemma} \label{lems2} Suppose the positive parameters $\epsilon_1,\epsilon_2$ and $\delta$ are chosen such that $\epsilon_2\le \epsilon_1|\log\epsilon_1|,\, \epsilon_1<\delta^2$ and  $\delta\ll 1$. Assume that $\ell^2(\mathbb{N})$-valued function $h$ is holomorphic in  $R_{I_{1+\delta},\epsilon_1}$, $h\in C(\overline{R_{I_{1+\delta},\epsilon_1}})$, $h\nequiv 0$,
\[
W:=\sup_{0<y<\epsilon_1}\int_{I_{1+\delta}} \|h(x+iy)\|^2dx, \quad \|h(k)\|^2\le  L (y^{-1}+(\epsilon_1-y)^{-1})\,, \quad k=x+iy\in R_{I_{1+\delta},\epsilon_1}\,,\quad L>2\,.
\] 
Then, we have
\begin{eqnarray*}
\inf_{0<y<\epsilon_2/2} \int_{I_{1-\delta}}\log \|h(x+iy)\|dx\ge 
(1+O(\epsilon_2\epsilon_1^{-1}))\left(\int_{I_1}\log\|h(x+i\epsilon_2)\|dx-\eta\right), \\
|\eta|<C\Bigl(\epsilon_2\epsilon_1^{-1}\Bigl(W^{0.5}+ |\log L| +|\log\epsilon_1|\Bigr)+(\delta W)^{0.5}\Bigr)\,.
\end{eqnarray*}
\end{lemma}
\begin{proof}
It is enough to prove 
\begin{equation}\label{p9}
 \int_{I_{1-\delta}}\log \|h(x)\|dx\ge 
(1+O(\epsilon_2\epsilon_1^{-1}))\left(\int_{I_1}\log\|h(x+i\epsilon_2)\|dx-\eta\right).
\end{equation}
Take $k=x+i\epsilon_2, x\in I_{1}$ and apply the mean-value inequality to the subharmonic function $\log \|h\|$ within $R_{I_{1+\delta},\epsilon_1}$.
We define
$
\Gamma_1=\{k:\Re k\in I_{1+\delta},\Im k=\epsilon_1\}\,,\,
\Gamma_2=\{k:\Im k\in (0,\epsilon_1), k\in \partial {R_{I_{1+\delta},\epsilon_1}}\}\,,\,
\Gamma_3=\{k:\Re k\in I_{1+\delta},\Im k=0\}\,.
$ Check Figure 4.

\vspace{1cm}

\begin{tikzpicture}
\draw (0,1) -- (13,1);
\draw (1,1) -- (1,4);
\draw (1,4) -- (12,4);
\draw (12,4) -- (12,1);

\draw (2,1) -- (2,4);

\draw (3,1) -- (3,1.5);
\draw (2,1.5) -- (11,1.5);
\draw (11,4) -- (11,1);

\draw (10,1.5) -- (10,1);

\draw (10.2,1.5)   node[anchor=north]  {$\epsilon_2$};

\draw (11.3,2.8)   node[anchor=north]  {$\epsilon_1$};

\draw (6,0.5)   node[anchor=north] 
 {Figure 4: $\epsilon_2\ll \epsilon_1< \delta^2\ll 1$};

\draw plot [smooth] coordinates {(1,1) (1,0.9) 
(1.4,0.8)
(1.5,0.7) 
(1.6,0.8)
 (2,0.9) (2,1) };
 
 \draw plot [smooth] coordinates {(2,1) (2,0.9) 
(2.4,0.8)
(2.5,0.7) 
 (2.6,0.8)
 (3,0.9) (3,1) };
 
\draw (1.5,0.6)   node[anchor=north]  {$\delta$};

\draw (2,0.8) node[anchor=north] {$-1$};

\draw (11,0.8) node[anchor=north] {$1$};

\draw (2.5,0.6)   node[anchor=north]  {$\delta$};
\draw (4.5,2)   node[anchor=north]  {$k$};
\draw[thick] (4.4,1.5)  circle (0.2mm);

\draw (6.5,0.8) node[anchor=north] {$0$};
\draw[thick] (6.5,1)  circle (0.2mm);

\end{tikzpicture}

We get
\begin{equation}\label{p1}
\int_{\Gamma_3}\log\|h\|d\omega_k\ge \log \|h(x+i\epsilon_2)\|-E_1-E_2\,,
\end{equation}
where 
\[
E_1=\int_{\Gamma_1}\log_+ \|h\|d\omega_k, \quad 
E_2=\int_{\Gamma_2}\log_+ \|h\|d\omega_k\,.
\]
One applies the given estimates on $h$ and the estimates on a  harmonic measure to bound $E_{1(2)}$:
\[
E_2\lesssim   (|\log L| +|\log\epsilon_1|)\delta^{-3}\epsilon^2_1\epsilon_2, \quad
E_1\lesssim \int_{I_{1+\delta}}\log_+\|h(\xi+i\epsilon_1)\| \frac{\epsilon_1^{-2}y}{\cosh (\pi \epsilon_1^{-1}(x-\xi))}  d\xi\,.
\]
Now, we integrate \eqref{p1} in $x$ over $I_{1}$ and recall that $\Gamma_3=I_{1+\delta}$.
That gives
\[
\int_{I_1} E_1dx\lesssim \int_{I_{1+\delta}}\log_+\|h(\xi+i\epsilon_1)\| \left(\int_{I_1}\frac{\epsilon_1^{-2}y}{\cosh (\pi \epsilon_1^{-1}(x-\xi))} dx\right) d\xi\le W^{\frac 12}y\epsilon_1^{-1}\,.
\]
Then,
\begin{eqnarray*}
\int_{I_{1}}\left(\int_{I_{1+\delta}}\log \|h\|d\omega_k \right)dx=\int_{I_{1+\delta}} \log\|h\|\left( \int_{I_{1}}\omega_k'dx\right)d\xi\le\\ (1+O(\epsilon_2\epsilon_1^{-1}))\int_{I_{1-\delta}} \log\|h\|dx+
O(\epsilon_2\epsilon_1^{-1}))\int_{I_{1-\delta}} \log_+\|h\|dx+
\int_{I_{1+\delta}\backslash I_{1-\delta}} \log_+\|h\|\left( \int_{I_{1}}\omega_k'dx\right)d\xi\le \\
 (1+O(\epsilon_2\epsilon_1^{-1}))\int_{I_{1-\delta}} \log\|h\|dx+C\epsilon_2\epsilon_1^{-1}W^{\frac 12}+C\int_{I_{1+\delta}\backslash I_{1-\delta}} \log_+\|h\|d\xi
\end{eqnarray*}
after we use the bound \eqref{new1} from Lemma \ref{new_l1}.
Finally, 
\[
\int_{I_{1+\delta}\backslash I_{1-\delta}} \log_+\|h\|d\xi\le C(\epsilon)W^{0.5}\delta^{0.5}
\]
by Cauchy-Schwarz inequality. Combining obtained estimates, we get the statement of the lemma.
\end{proof}
\bigskip

\section{Appendix 2: rough bounds on Green's function}

We need the following  standard bounds ``a la Combes-Thomas''  (see, e.g., \cite{ct}) for Green's function $G(r,\rho,k^2)$ of $H=H^{(0)}+V$. In this section, we assume that $I$ is a fixed closed interval in $\R^+$ and $k\in R_{I,1}$.
%\begin{lemma}Suppose $\|V\|_{L^\infty(\R^+)}<\infty$ and $k\in \Pi_1$. Then, there is $c=c(V,I)>0$ such that
%\[
%\|G(x,y,k^2)\|<C(V,I)\frac{e^{-c|x-y|\Im k}}{\Im k}\,.
%\]
%\end{lemma}
%\begin{proof}
%This is Combes-Thomas estimate. 
%\end{proof}
\begin{lemma}Suppose $\|V\|_{L^\infty(\R^+)}<\infty$. Then, we have
\begin{equation}\label{g1}
\|G(r,\rho,k^2)\|\le C'_Ie^{-0.5(\Im k) |r-\rho|}
\end{equation}
for all $k\in R_{I,1}, \Im k>C_I\|V\|_{L^\infty(\R^+)}$ with some $C_I>0$ and $C_I'>0$.
\end{lemma}
\begin{proof}This is immediate from the analysis of perturbation identity for the Green's kernel $G$:
\[
G(r,\rho,k^2)=G^{(0)}(r,\rho,k^2)-\int_0^\infty  G^{(0)}(r,\xi,k^2)V(\xi)G(\xi,\rho,k^2)d\xi\,.
\]
Multiply the both sides by $e^{0.5(\Im k)|r-\rho|}$ and apply the contraction mapping principle in $L^\infty(\R^+\times \R^+)$. We use \eqref{fi} to get
\[
e^{0.5(\Im k) |r-\rho|}\int_0^\infty e^{-(\Im k) |r-\xi|}\|V(\xi)\|e^{-0.5(\Im k) |\xi-\rho|}d\xi\le 4\|V\|_{L^\infty(\R^+)}(\Im k)^{-1}
\]
and \eqref{g1} follows provided  $\Im k>C_I\|V\|_{L^\infty(\R^+)}$ with suitable $C_I$.
\end{proof}

Finally, we can focus on the lemma we need in the main text.

\begin{lemma} \label{l63} Let $\|V\|\le \lambda (1+r)^{-\gamma},\, H=H^{(0)}+V,\, k\in R_{I,1}$, where $I$ is a closed interval in $\R^+$, $\gamma\in (0,1)$, and $T>1$. Then, there are positive $T$-independent constants $C,C_1$ and $c$ such that
\[
\|G(r,\rho,k)\|<Ce^{-c(\Im k) |r-\rho|}
\]
for $ \Im k>C_1T^{-\gamma}$, $0.5T<r<T$, and $0.5T<\rho<T$. 
\end{lemma}
\begin{proof}Define $H'=-\partial^2_{rr}+V\cdot \chi_{r>\frac 14 T}$.  By the previous lemma, the corresponding Green's kernel $G'$ satisfies the bound
\begin{equation}\label{h1}
\|G'(r,\rho,k)\|\le Ce^{-0.5(\Im k)|r-\rho|}
\end{equation}
if $\Im k>C_1T^{-\gamma}$. Next, we again write the second resolvent identity
\[
G(r,\rho,k^2)=G'(r,\rho,k^2)-\int_0^{\frac 14  T} G(r,\xi,k^2)V(\xi)G'(\xi,\rho,k^2)d\xi\,.
\]
For the first term, we use \eqref{h1}. To estimate the second one, we apply a general bound: for every $h\in \mathcal{H}$, one has
$
\|R_{k^2}h\|_{L^\infty(\R^+)}\le C(\|R_{k^2}h\|_{L^2(\R^+)}+\|(R_{k^2}h)''\|_{L^2(\R^+)})\le C_{I,\lambda}(\Im k)^{-1}\|h\|_{L^2(\R^+)}
$
which follows from Sobolev's embedding, the equation for $R_{k^2}h$, and the Spectral Theorem. Then, since $r,\rho\in [0.5T,T]$, one deduces
\[
\left\|\int_0^{\frac 14  T} G(r,\xi,k^2)V(\xi)G'(\xi,\rho,k^2)d\xi\right\|\le C_{I,\lambda}(\Im k)^{-1}\left(\int_0^{\frac 14  T}e^{-(\Im k)|\xi-\rho|}d\xi\right)^{\frac 12}\le C_{I,\lambda}(\Im k)^{-2} e^{-0.1 (\Im k)T}\,.
\]
Since $\Im k>C_1T^{-\gamma}$ and $\gamma\in (0, 1)$, we have $(\Im k)^{-2} e^{-0.1(\Im k)T}<Ce^{-c_1(\Im k)T}$ with positive $c_1$. The result now follows because $e^{-c_1(\Im k)T}\le e^{-c(\Im k)|\rho-r|}$ with positive $c$ provided that $0.5T<r,\rho<T$.
\end{proof}

\end{document}